\documentclass[a4paper,12pt]{article}
\usepackage[T2A]{fontenc}
\usepackage[english,russian]{babel}
\usepackage{amsfonts,amssymb,amscd,amsmath}
\usepackage[dvips]{graphicx}
\usepackage[12pt]{extsizes}
\usepackage[right=25mm,
left=25mm,top=30mm,
bottom=40mm]{geometry}

\begin{document}

\begin{center}{Resonance vibrations of impact oscillator 
with biharmonic 
excitation}\end{center}

\vskip 12pt\noindent

\centerline{\bf Vladimir~Burd}

\vskip 12pt\noindent

{Department of Mathematics, Yaroslavl State
University, Russia}

\bigskip

Abstract.
We consider a damped impact oscillator subject
to the action of a biharmonic force. The conditions 
for the existence and stability of almost periodic
resonance solutions are investigated.



\section{Introduction}
Representations of vibro-impact processes 
in terms of integrals of motion for 
vibro-impact system with one 
degree of freedom is determined by the 
physics 
of vibro-impact system. The corresponding 
variables (impulse - phase) are very natural 
for the  vibrating
processes accompanied by impacts. The 
impulse can be interpreted as the force 
characteristics of impact. 
The phase is the instant of 
collision.
\par
In case of an instantaneous impact the 
representations of 
periodic vibro-impact processes can be 
written via so-called periodic Green's functions.
Periodic Green's functions are determined only by 
the linear part of the system. They appear 
below or 
in the form of Fourier series or finite 
relations 
on intervals of periodicity. 
The impact force is
represented by the Dirac $\delta$-function 
and the values depend on the impulse and 
phase of the impact.  
The periodic Green's functions are calculated 
using well known methods [2]. Similar 
approach was proposed by V.I. Babitsky 
and M.Z. Kolovsky [3].
\par
We want to give a correct description of resonant 
modes in vibro-impact system with one degree of 
freedom.
\par
We chose to study a model with elastic impact and
viscous damping. Such a model 
is discussed in [2], [14]. 
\par
A model that takes into account the energy
loss at impact is proposed in [2]. It is not clear 
how to explore this model on a 
mathematical level 
of rigor.
\par
Our study is based on a 
combination of the method
of averaging on infinite interval [6] and 
the methods of singular perturbation theory 
[10,13].
\par
We will study the motion of 
conservative vibro-impact 
systems with one degree of freedom 
and subject to small biharmonic 
perturbations.
This problem is reduced
to study two-dimensional 
 systems with fast rotating
phase and slow varying coefficients.
\par
Periodic perturbations of smooth two-dimensional 
systems with fast rotating phase and slow varying
coefficients were studied in [9]. 
The conditions
of closeness of solutions of exact and averaged
equations on a finite asymptotically large time
interval are established. Resonance almost periodic
oscillations in such systems are investigated in 
[5] and [6].
\par
Recently, much attention has been focused on 
studies of dynamics of 
nonlinear systems perturbed 
by a biharmonic external force
with different 
frequencies (see, for example [1, 4, 7, 8]).

\section{Conservative impact oscillator}

Here we follow the method from the book [2].
Consider a linear oscillator 
$$
\ddot x+\Omega^2 x=0.
$$
\par
At the point $x=\Delta$, we arrange 
an immovable limiter and assume
that once the coordinate  $x$ 
reaches the value
$\Delta$, an instant elastic impact 
occurs in the system
so that if $x=\Delta$ at the time 
instant  $t_\alpha$,
then the relation
\begin{equation}
\dot x(t_\alpha-0)=-\dot x( t_\alpha+0) 
\label{y}
\end{equation}
holds.
\par
If $\Delta>0$ and the energy level in the linear 
system is insufficient to attain 
the level  $x=\Delta$, then linear oscillations 
with the frequency $\Omega$ take place. In the 
presence of collisions, the oscillation 
frequency satisfies the inequality $\omega>\Omega$ and
rises as the energy rises but does not exceed
the value  $2\Omega$. Hence
\begin{equation}
\Omega<\omega<2\Omega,\quad
\Delta>0. 
\label{(a)}
\end{equation}
\par
If $\Delta<0$, the 
oscillation frequency $\omega$ obeys the
inequality
\begin{equation}
2\Omega<\omega<\infty,\quad
\Delta<0.     
\label{b}
\end{equation}
\par
At $\Delta=0$ the image point 
passes any phase
trajectory for the same time with 
the doubled velocity $2\Omega$ so that
\begin{equation}
\omega=2\Omega,\quad \Delta=0.
\label{c}
\end{equation}
\par
Condition (1) suggests that 
variation of the impulse $\Phi_0$
in the neighborhood of the impact 
instant  $t_\alpha$ takes the
form
$$
J=\dot x_--\dot x_+=2\dot x_-,
\quad \dot x_->0,
$$
where $\dot x_\alpha=\dot x(t_\alpha \mp 0)$.
\par
The resulting force becomes localized 
at $t=t_\alpha$. Hence
\begin{equation}
\Phi_0\left|_{t=t_\alpha}\right. =J
\delta(t-t_\alpha) 
\label{d}
\end{equation}
and
$$
\int\limits_{t_\alpha-0}^{t_\alpha
+0}\Phi_0dt=J.
$$
Impacts occur periodically when 
$t_\alpha=t_0+\alpha T$, where
$\alpha$ is an integer, and $T$ is 
the period between impacts
calculated by the equality $T=2\pi
\omega^{-1}$ and (2)--(4). Thus, for 
$-\infty<t<\infty$, we obtain a $T$-periodic
continuation of (5)
$$
\Phi_0=J\delta_T(t-t_0),
$$
where $\delta_T(t)$ is the $T$-periodic 
$\delta$-function.
\par
Solution of the equation
\begin{equation}
\ddot x+\Omega^2x+
\Phi_0(x,\dot x)=0  
\label{e}
\end{equation}
is understood as the $T$-periodic 
function $x(t)$ such that its
substitution into this equation 
transforms it into a correct
equality (from the viewpoint of the 
theory of distributions) of the form
$$
\ddot x+\Omega^2x+J\delta_T(t-t_0)=0,
$$
where $t_0$ is an arbitrary 
constant, and for all $\alpha= 0,\pm 1,\dots$
$$
x(t_0+\alpha T)=\Delta,\quad
J=2\dot x_-(t_0+\alpha T).
$$
At the same time, the restrictions
$$
x(t)\le\Delta,\quad \dot x_->0 
$$
are fulfilled, and the periods of 
oscillations, depending on the
sign of $\Delta$, fit the frequency 
ranges of (2)--(4).
\par
To describe the solution analytically, 
we assume $t_0=0$. In this
case, for $0\le t<T_0$, the solution 
of equation (6) has the form
$$
\begin{array}{l}
x(t)=-J\kappa[\omega_0(J)(t-t_0),
\omega_0(J)],\quad
\kappa(t,\omega_0)=\frac{1}
{2\Omega}\frac{\cos[\Omega(t-T_0/2)]}
{\sin(\Omega T_0/2)}=\\\frac{\omega_0}
{2\pi\Omega^2}+\frac{\omega_0}{
\pi}\sum_{k=1}^\infty\frac{
\cos k\omega_0 t}{\Omega^2-k^2
\omega_0^2},\quad
J(\omega_0)=-2\Omega\Delta\tan
\frac{\Omega T_0}{2},\quad J\ge 0,
\end{array}
$$
and the third relation here 
determines the smooth dependence
$\omega_0(J)$ at $\Delta\ne 0$, whereas $\omega_0=2\Omega$ at
$\Delta=0$.
\par
Geometric conditions of an impact 
result in frequency intervals (2)--(4). Note 
that when $\Delta=0$, the solution
$x(t)$ for $0\le t<\pi/\Omega$ takes the form
$$
x(t)=-\frac{J}{2\Omega}
\sin\Omega t,
$$
where  $J$ is a frequency-independent arbitrary 
constant.

\section{Perturbed impact oscillator}

Now consider a perturbed impact oscillator
\begin{equation}
\ddot x+\Omega^2 x+\Phi_0(x,\dot x)=\varepsilon[f(t,\tau)-\gamma\dot x],
\label{1}
\end{equation}
where $\varepsilon>0$ is a small parameter,
$\tau=\varepsilon t$ is a slow time, 
$\gamma>0$ is a constant. The function 
$\Phi_0(x,\dot x)$ describes the force of impact.
Representation of this function is given in the
previous section.
\par
We consider the perturbations 
of two types. First
a biharmonic perturbation is 
the sum of two small
periodic forces with close 
frequencies.
The corresponding perturbation 
has the form
$$
f(t,\tau)=a_1\sin\nu t+a_2\sin(\nu t+
\Gamma\tau),
$$ 
where $a_1,~a_2,~\nu,~\Gamma$ are real positive
numbers. The perturbation function $f(t,\tau)$ 
can be represented as
\begin{equation}
f(t,\tau)=E(\tau)\sin(\nu t+
\beta(\tau)),\quad 
\tau=\varepsilon t, 
\label{g}
\end{equation}
where
$$
\begin{array}{l}
E(\tau)=\sqrt{a_1^2+2a_1a_2\cos\Gamma
\tau+a_2^2},\quad \cos\beta(\tau)=
\frac{a_1+a_2\cos\Gamma\tau}{E(\tau)},
\\
\sin\beta(\tau)=\frac{a_2
\sin\Gamma\tau}{E(\tau)}.
\end{array}
$$
The function (8) is periodic in $t$ with the period
$2\pi/\nu$ and periodic in $\tau$ with the period 
$2\pi/\Gamma$. The function $E(\tau)$ is strictly 
positive if $a_1\ne a_2$, which will be assumed. 
Consequently
the function $f(t,\varepsilon t)$ is almost periodic
function (see [6]) with two basic frequencies.
\par
Second perturbation is a biharmonic force with very
different frequencies. This force has the form
\begin{equation}
f(t,\tau)=A\sin(\nu t+\theta)+B\sin\Gamma\tau,
\quad \tau=\varepsilon t.
\label{h}
\end{equation}
where $A$,~$B$,~$\nu$,~$\Gamma$,~$\theta$ are 
real numbers.
\par
We assume  
$\psi=\omega_0(J)t$ and transform
equation (7) into the system in 
terms of the variables $J,\psi$
(impulse-phase), by making a change
\begin{equation}
\begin{array}{l}
x=-J\kappa[\psi,\omega_0(J)],\\
\dot x=-J\omega_0(J)\kappa_\psi
[\psi,\omega_0(J)],\end{array}
\label{2}
\end{equation}
where
$$
\kappa[\psi,\omega_0(J)]=\kappa(\psi,J)=
\omega_0(J)^
{-1}\left[\frac{1}{2\pi\Omega_0^2}
+\frac{1}{\pi}\sum_{k=1}^\infty
\frac{\cos k\psi}{\Omega_0^2-k^2}
\right],\quad \Omega_0=\frac{\Omega}
{\omega_0(J)}.
$$
It follows from the theory of Fourier series 
(see, for example, [12, chapter 4]) that the
function $\kappa_\psi[\psi,\omega_0(J)]$ has a 
finite discontinuities at the points 
$\psi=2l\pi$, where $l$ is an integer, and 
continuous at all other points. Moreover the
function 
$\kappa_\psi[\psi,\omega_0(J)]$ 
has a derative 
at interior 
points of the
intervals $[2\pi l,2\pi(l+1)]$.
Fourier series of
function $\kappa_\psi[\psi,\omega_0(J)]$ has a 
form
$$
-\omega_0(J_{pq})^{-1}\frac{1}{\pi}
\sum_{k=1}^\infty\frac{k\sin k\psi}
{\Omega_0^2-k^2}.
$$
\par
Therefore change  (10) is not smooth 
at $\psi=2l\pi$.
Thus, in new variables, impacts occur  
when  $\psi=2l\pi$. Making the
substitution (10), we arrive at 
the system 
\begin{equation}
\begin{array}{l}
\frac{dJ}{dt}=-4\varepsilon
\omega_0(J)[f(t,\tau)+\gamma J\omega_0(J)
\kappa_\psi(\psi,J)]\kappa_\psi(\psi,J),
\\
\frac{d\psi}{dt}=\omega_0(J)-
4\varepsilon
\omega_0(J)J^{-1}[f(t,\tau)+J\omega_0(J)\kappa_\psi
(\psi,J)](-J\kappa(\psi,J))_J.
\end{array}
\label{3}
\end{equation}
A detailed derivation of the system (11) is
contained in the book [2, chapter 2].
This is a system with a fast 
rotating phase, where the right-hand
sides are periodic in  $\psi$ and 
have finite discontinuities
at the points $\psi=2l\pi$. The dependencies
$\omega_0(J)$ have the form
$$
\begin{array}{l}
\omega_0(J)=\frac{\pi\Omega}{\pi-
\arctan[J/(2\Omega\Delta)]},
\quad \Delta>0,\quad \Omega<
\omega_0<2\Omega,\\
\omega_0(J)=-\frac{\pi\Omega}
{\arctan[J/(2\Omega\Delta)]},
\quad \Delta<0,\quad 2\Omega<
\omega_0<\infty,\\
\omega_0=2\Omega=const,\quad
\Delta=0.\end{array}
$$
System (11) is a system with two 
slow variables $J,\tau$ and two fast variables 
$\psi,t$. The existence and stability of
stationary resonance  solutions in such systems
will be investigated. 
\section{Construction of averaged equations}

We will use the method
of averaging on an infinite interval (see [6]). 
The major part in all the problems related to 
the principle of averaging is in 
the changes of variables. The
changes allow to eliminate fast variables 
from equations of motion within the given accuracy 
and thus separate slow motion from the fast one.
\par
Let $J_{pq}$ be a solution of the equation
$$
\omega_0(J_{pq})=\frac{q}{p}\nu.
$$
where $p$,~$q$ are relatively prime integers.
By making a change $\psi=\varphi+\frac{q}{p}\nu t$
we transform system (11) into
\begin{equation}
\begin{array}{l}
\frac{dJ}{dt}=-4\varepsilon
\omega_0(J)[f(t,\tau)+\gamma J\omega_0(J)
\kappa_\psi(\varphi+\frac{q}{p}\nu t,J)]\kappa_\psi(
\varphi+\frac{q}{p}\nu t,J),
\\
\frac{d\varphi}{dt}=\omega_0(J)-\frac{q}{p}\nu -
4\varepsilon
\omega_0(J)J^{-1}[f(t,\tau)+\gamma J\omega_0(J)\kappa_\psi
(\varphi+\frac{q}{p}\nu t,J)]\times\\
(-J\kappa(\varphi+\frac{q}{p}\nu t,J))_J.
\end{array}
\label{4}
\end{equation}
The point $J=J_{pq}$ is resonant point in system
(12).
\par
Assume that the resonance is non-degenerate, i.e.,
\begin{equation}
\frac{d\omega_0}{dJ}\biggl|_{J=J_{pq}}=\omega_0'(J_{pq})
\ne 0.
\label{k} 
\end{equation}
We shall study the behavior of solutions of 
system (12) in the 
$\mu=\sqrt\varepsilon$-neighborhood 
of the resonance point $J_{pq}$. We make a change
$$
J=J_{pq}+\mu z
$$
and expand the right-hand side of 
system (12) in terms of the powers of $\mu$.
As a result, we obtain the system
\begin{equation}
\begin{array}{l}
\frac{dz}{dt}=\mu F_0(t,\tau,\varphi,J_{pq})+
\mu^2F_1(t,\tau,\varphi,J_{pq})z+O(\mu^3),\\
\frac{d\varphi}{dt}=\mu\omega_0'(J_{pq})z+\frac{1}{2}
\mu^2\omega_0''(J_{pq})z^2+\mu^2G_0(t,\tau,\varphi,J_{pq})
+O(\mu^3)
\end{array}
\label{5}
\end{equation}
where 
$$
\begin{array}{l}
F_0(t,\tau,\varphi,J_{pq})=-4\omega_0(J_{pq})[f(t,\tau)+
\gamma J_{pq}\omega_0(J_{pq})
\kappa_\psi(\varphi+\frac{q}{p}\nu t,J+{pq})]\times\\
\kappa_\psi(\varphi+\frac{q}{p}\nu t,J_{pq}),
\end{array}
$$
$$
\begin{array}{l}
F_1(t,\tau,\varphi,J_{pq})=-4\frac{d}{dJ}\left\{
\omega_0(J)[f(t,\tau)+\gamma J\omega_0(J)
\kappa_\psi(\varphi+\frac{q}{p}\nu t,J)]\times\right.\\
\left.\kappa_\psi(
\varphi+\frac{q}{p}\nu t,J)\right\}\biggl|_
{J=J{pq}}. 
\end{array}                                           
$$
$$
\begin{array}{l}
G_0(t,\tau,\varphi,J_{pq})=-4\omega_0(J_{pq})
J^{-1}_{pq}[f(t,\tau)+\gamma J_{pq}\omega_0(J_{pq})
\kappa_\psi(\varphi+\frac{q}{p}\nu t,J_{pq})]\times\\
(-\kappa(\varphi+\frac{q}{p}\nu t,J_{pq})-J_{pq}
\kappa_J(\varphi+\frac{q}{p}\nu t,J_{pq}).
\end{array}
$$
System (14) contains only one 
fast variable $t$. We now make
the standard change of the method 
of  averaging in order to eliminate 
the fast variable in the right-hand 
side of system (14) up to the 
accuracy of the terms of order 
$\mu^2$. This change is sought in the 
form
$$
z=\xi+\mu u_1(\eta,t,\tau)+
\mu^2u_2(\eta,t,\tau)\xi,\quad
\varphi=\eta+\mu^2v_2(\eta,t,
\tau),
$$
where the functions  $u_i(\eta,t,
\tau),\,(i=1,2),\,v_2(\eta,t,\tau)$ are 
periodic in $t$,~$\tau$ with period $2\pi/\nu$ and
$2\pi/\Gamma$ accordingly. 
\par
The change results in the system
\begin{equation}
\begin{array}{l}
\frac{d\xi}{dt}=\mu f_0(\eta,\tau)
+\mu^2f_1(\eta,\tau)\xi+
O(\mu^3),\\
\frac{d\eta}{dt}=\mu\omega_0'(J_{pq})\xi+\frac{1}{2}
\mu^2\omega_0''(J_{pq}))\xi^2+\mu^2
g_0(\eta,\tau)+O(\mu^3),
\end{array}
\label{6}
\end{equation}
where $f_0(\eta,\tau)$,~$f_1(\eta,\tau)$,~$g_0(\eta,\tau)$ 
are defined as the mean values over $t$:
$$
\begin{array}{l}
f_0(\eta,\tau)=\frac{\nu}{2\pi}\int\limits_0^
{2\pi/\nu}F_0(t,\tau,\eta,J_{pq})dt, \\
f_1(\eta,\tau)=\frac{\nu}{2\pi}\int\limits_0^{2\pi/\nu}
F_1(t,\tau,\eta,J_{pq})dt,\\
g_0(\eta,\tau)=\frac{\nu}{2\pi}\int\limits_0^
{2\pi/\nu}G_0(t,\tau,\eta,J_{pq})dt.
\end{array}
$$
The functions $u_i(\eta,t,\tau),~(i=1,2)$,~
$v_2(\eta,t,\tau)$ are defined as periodic solutions
in $t$ with the zero mean value from the equations
$$
\begin{array}{l}
\frac{\partial u_1}{\partial t}=F_0(t,\tau,\eta,J_{pq})
-f_0(\eta,\tau),
\\
\frac{\partial u_2}{\partial t}=F_1(t,\tau,\eta,J_{pq})
-u_{1\eta}(\eta,t,\tau)\omega_0(J_{pq})-f_1(\eta,\tau),
\\
\frac{\partial v_2}{\partial t}=G_0(t,\tau,\eta,J_{pq})-
-u_1(\eta,t,\tau)\omega_0(J_{pq})-g_0(\eta,\tau).
\end{array}
$$
\par
Functions $f_0(\eta,\tau)$,~$f_1(\eta,\tau)$,~
$g_0(\eta,\tau)$ are smooth in $\eta$ for considered
perturbations. The function $f_0(\eta,\tau)$
is mean value over t of the function 
$F_0(t,\tau,\varphi,J_pq)$. 
After averaging only a single term remains from
the infinite sum.
Calculation of $f_0(\eta,\tau)$ will be 
demonstrated below (see (34)).
 Similar assertions hold
for functions  $f_1(\eta,\tau)$,~$g_0(\eta,\tau)$.
\par 
System (15) at the time $\tau$ is a singularly 
perturbed system in the following form
\begin{equation}
\begin{array}{l}
\mu \frac{d\xi}{d\tau}=f_0(\eta,\tau)+\mu
f_1(\eta,\tau)\xi+O(\mu^2),\\
\mu\frac{d\eta}{d\tau}=\omega_0'(J_{pq})\xi+
\frac{1}{2}\mu\omega_0''(J_{pq}))\xi^2+\mu
g_0(\eta,\tau)+O(\mu^2),
\end{array}
\label{7}
\end{equation}
\section{Existence and stability of almost periodic
solutions}
\label{w}

Let there exist periodic function $\eta_0(\tau)$ 
such that
\begin{equation}
f_0(\eta_0(\tau),\tau)\equiv 0,\quad 0< 
\eta_0(\tau)<2\pi.
\label{8}
\end{equation}
In this case, the degenerate system derived from 
(16) at $\mu=0$ has the solution
\begin{equation}
\xi=0,\quad \eta=\eta_0(\tau).
\label{9}
\end{equation}
Linearizing the right-hand side of 
system (16) at $\mu=0$ on
solution (18) yields the matrix
$$
A_0(\tau)=\left(\begin{array}{cc}
0&f_{0\eta}(\eta_0,\tau)\\\omega_0'
(J_{pq})&0\end{array}\right).
$$
If
\begin{equation}
\omega_0'(J_{pq})f_{0\eta}
(\eta_0,\tau)>\sigma_0>0,\quad \sigma_0=const,
\quad \tau\in(-\infty,\infty),
\label{10}
\end{equation}
then the matrix $A_0(\tau)$ has real eigenvalues of 
different signs. In this case, as we know 
[6, chapter 8], in the system
\begin{equation}
\mu\frac{du}{d\tau}=A_0(\tau)u
\label{l}
\end{equation}
for sufficiently small $\mu$, the space of solutions
$U(\mu)$ can be represented in the form
$$
U(\mu)=U_+(\mu)+U_-(\mu).
$$
For the solutions $u_+(\tau,\mu)\in U_+(\mu)$ the
inequality 
$$
|u_+(\tau,\mu)|\le M_+\exp\left[-\frac{\gamma_+}{\mu}
(\tau-s)\right]|u_+(s,\mu)|,\quad (-\infty<s<\tau<\infty)
$$
holds, and for $u_-(\tau,\mu)\in U_-(\mu)$ the 
following inequality holds:
$$
|u_-(\tau,\mu)|\le M_-\exp\left[-\frac{\gamma_-}{\mu}
(\tau-s)\right]|u_+(s,\mu)|,\quad (-\infty<\tau<s<\infty).
$$
Here $M_+,M_-,\gamma_+,\gamma_-$ are positive 
constants and $|\cdot|$ is a norm in $\mathbb{R}^2$.
It follows from an estimation of the solutions of 
system (20) that the solution of this system is
unstable for sufficiently small $\mu$ if the space
$X_-(\mu)$ of the initial conditions of the solutions 
from $U_-(\mu)$ is non-trivial.
Hence, an inhomogeneous system
$$
L(\mu)=\mu\frac{dz}{d\tau}-A_0(\tau)z=
f(\tau),
$$
where $f(\tau)$ is a periodic two-dimensional 
function, for sufficiently small $\mu$
has a unique periodic solution. This solution is
represented as 
$$
z(\tau,\mu)=L^{-1}(\mu)f=\frac{1}{\mu}
\int\limits_{-\infty}
^\infty K(\tau,s,\mu)f(s)ds,
$$
where
\begin{equation}
|K(\tau,s,\mu)|\le M\exp\left(-\frac
{\gamma}{\mu}|\tau-s|\right)\quad
(-\infty<\tau,s<\infty),
\label{m}
\end{equation}
and $M$,~$\gamma$ are positive constants.
\par
We transform system (16) using a change
$$
u=\eta-\eta_0(\tau)
$$
and  write the obtained system in the vector form 
$(z=(\xi,u))$ 
\begin{equation}
\mu\frac{dz}{d\tau}=A_0(\tau)z+
H(z,t,\tau,\mu).
\label{n}
\end{equation}
The components $H(z,t,\tau,\mu)$ 
have the form
$$
\begin{array}{l}
f_0(u+\eta_0,\tau)-f_{0\eta}(\eta_0,\tau)u+
\mu f_1(u+\eta_0,\tau)\xi+O(\mu^2),
\\
-\mu\frac{d\eta_0}{d\tau}+
\frac{1}{2}\mu\omega_{xx}(x_0,
\tau)\xi^2+\mu g_0(u+\eta_0,\tau)+O(\mu^2).
\end{array}
$$
Evidently, the following inequality is valid
\begin{equation}
|H(0,t,\tau,\mu)|\le p(\mu),
\label{o}
\end{equation}
where $p(\mu)\to 0$ as $\mu\to 0$. 
The components of the vector function $H(z,t,
\tau,\mu)$ are differentiable at z in a sufficiently
small neighborhood of the point $(0,0)$.
Therefore the following inequality holds
\begin{equation}
|H(z_1,t,\tau,\mu)-H(z_2,t,\tau,\mu)|
\le p_1(r,\mu)|z_1-z_2|, 
\label{p}
\end{equation}
where $p_1(r,\mu)\to 0$ as $r\to 0$. 
\par
The vector function $H(z,t,\varepsilon t,\mu)$ is
an almost periodic function of variable $t$.
The problem of almost periodic 
solutions of system (22) is
equivalent to the problem of 
solvability of the system of integral equations
\begin{equation}
z(\tau,\mu)=\frac{1}{\mu}
\int\limits_{-\infty}^\infty
K(\tau,s,\mu)H(z,s,\mu)ds=\Pi(z,\mu).
\label{q}
\end{equation}
\par
From inequalities (23), (24) follow that
the successive approximations (see, for example, 
[11])
$$
z_0(\tau,\mu)=\Pi(z,0),\quad z_j(\tau,\mu)=
\Pi(z_{j-1}(\tau,\mu),\mu),\quad j=1,2,\dots
$$
for sufficiently small $\mu$ converge uniformly
on interval $(-\infty,\infty)$ to a unique almost
periodic solution $z^*(\tau,\mu)$ of the system of
integral equations (25).
The solution $z^*(\tau,\mu)$ tends to  $(0,0)$ 
uniformly with respect to  $t$ as $\mu\to 0$.
Hence, system (22),
for sufficiently small $\mu$, has a 
unique almost periodic solution
$z_*(\tau,\mu)$. In its turn, 
system (16), for sufficiently small $\mu$, has a 
unique almost periodic solution.
Therefore, system (15), for sufficiently 
small $\mu$, has a unique almost periodic 
solution.
To investigate the stability of the 
almost periodic solution
$z_*(\tau,\mu)$ of system (22), we make a change 
$z=z_*(\tau,\mu)+y(\tau,\mu)$ and 
obtain the system
\begin{equation}
\mu\frac{dy}{d\tau}=A_0(\tau)y+
H_1(y,t,\tau,\mu),
\label{r}
\end{equation}
where
$$
H_1(y,\psi,\tau,\mu)=H(z_*
(\tau,\mu)+y(\tau,\mu),t,
\tau,\mu)-H(z_*(\tau,\mu),
\psi,\tau,\mu).
$$
The problem of the stability of the 
almost periodic solution $z_*(\tau,\mu)$ is 
reduced to the problem of the stability of the 
zero solution of system  (26). 
Taking into account the exponential estimates 
on the solutions of system (20), based on the 
Liapunov's theorem of the instability by first 
approximation, we obtain that, for sufficiently
small $\mu$, the zero solution of 
system (16) is unstable. Hence, the almost periodic 
solution $z_*(\tau,\mu)$ of system (22) is unstable.
We shall state this result as a theorem 
applied to system (14). 

\bigskip

{\it
Let 
$J_{pq}$ be a constant such that
$$
\omega_0(J_{pq})=0,
\quad
\omega_0'(J_{pq})\ne 0.
$$
Let there exists a periodic 
function $\eta_0(\tau)$ ($0<\eta_0(\tau)<2\pi$)
such that
$$
f_0(\eta_0(\tau),\tau)=0
$$
and inequality (19) holds:
$$
\omega_0'(J_{pq})f_{0\eta}
(\eta_0,\tau)>\sigma_0>0,\quad \sigma_0=const,
\quad \tau\in(-\infty,\infty).
$$
In this case, system (14), for sufficiently small 
$\mu$, has a unique unstable
almost periodic solution.}

\bigskip                 

Hence, given the  
conditions of Theorem 1 are met in the
$\mu$-neighborhood of the resonance 
point $J_{pq}$, there exists
a unique unstable almost periodic 
solution.
\par
Now let, instead of inequality (19), 
the contrary inequality holds
\begin{equation}
\omega_0(J_{pq})f_{0\eta}
(\eta_0,\tau)<\sigma_1<0, \quad \sigma_1=const,
\quad -\infty<\tau<\infty.
\label{12}
\end{equation}
If inequality (27) holds, the 
eigenvalues of the matrix
$A_0(\tau)$ for all $\tau$ are purely 
imaginary. Now we need to
consider averaged equations of 
higher approximations (we used only 
the first approximation in satisfying 
inequality (19)).
\par 
We consider a narrower neighborhood of
the resonance point $J_{pq}$. We make a change
in system (15)
$$
\begin{array}{l}
\xi=\mu u(t)+\mu\xi_0(\tau)+\mu^3u_3(t,\tau), \\
\eta=\eta_0(\tau)+\mu v(t)+\mu^2v_0(\tau),
\end{array}
$$
where $\eta_0(\tau)$ satisfies equation (17),
periodic functions $\xi_0(\tau)$ and $v_0(\tau)$
are solutions of the equations
$$
\begin{array}{l}
\frac{d\eta_0}{d\tau}=\omega_0'(J_{pq})\xi_0(\tau)+
g_0(\eta_0(\tau),\tau),\\
\frac{d\xi_0}{d\tau}=f_{0\eta}(\eta_0(\tau),\tau)v_0
(\tau)+f_1(\eta_0(\tau),\tau)\xi_0(\tau)+\\\langle
F_1(t,\tau,\eta_0(\tau),J_{pq})u_1(\eta_0(\tau),t,\tau)\rangle_t
+\langle F_{0\varphi}(t,\tau,\eta_0(\tau),J_{pq})
v_2(\eta_0(\tau),t,\tau)\rangle_t
\end{array}
$$
respectively. Here $\langle\cdot\rangle_t$ is 
mean value in $t$ for fixed $\tau$.
These equations can be solved using
the inequalities (19) and (27). The functions
$u_3(t,\tau)$ is defined as periodic solution 
with zero mean value from the equation
$$
\begin{array}{l}
\frac{\partial u_3}{\partial t}=F_{0\varphi}(t,\tau,
\eta_0(\tau),J_{pq})v_2(\eta_0(\tau),t,\tau)+
F_1(t,\tau,\eta_0(\tau),J_{pq})u_1(\eta_0(\tau),t,
\tau)-\\
\langle F_{0\varphi}
(t,\tau,\eta_0(\tau),J_{pq})v_2(\eta_0(\tau),
t,\tau)\rangle_t-\langle F_1(t,\tau,\eta_0(\tau),J_{pq})
u_1(\eta_0(\tau),t,\tau)\rangle_t.
\end{array}
$$
The change transforms the system into
\begin{equation}
\begin{array}{l}
\frac{du}{dt}=\mu a(\tau)v+\mu^2
[b(\tau)u+e(\tau)v^2]+O(\mu^3),
\\
\frac{dv}{dt}=\mu c(\tau)u+
\mu^2d(\tau)v+O(\mu^3),
\end{array}
\label{13}
\end{equation}
where
\begin{equation}
\begin{array}{l}
a(\tau)=f_{0\eta}(\eta_0(\tau),\tau),\quad
b(\tau)=f_1(\eta_0(\tau),\tau), \quad
c(\tau)=\frac{1}{2}f_{0\eta\eta}(\eta_0(\tau),\tau),\\
d(\tau)=\omega'_0(J_{pq}),\quad
e(\tau)=g_{0\eta}(\eta_0(\tau),\tau).
\end{array}
\label{14}
\end{equation}
Condition (27) in new notation takes the form
$$
a(\tau)c(\tau)<\sigma_1<0.
$$
As was noted, it follows from this condition that 
the eigenvalues of the 
matrix of the first approximation
are purely imaginary for all $\tau$. We reduce 
system (28) to "standard form", i.e., to the 
form where the matrix of the first approximation 
is zero. 
A detailed description of reduction of system
(28) is contained in [6, chapter 16].
\par
We obtain the system
\begin{equation}
\begin{array}{l}
\frac{dA}{d\tau}=\frac{1}{2}
\left[b(\tau)+e(\tau)-
\frac{h'(\tau)}{h(\tau)}
\right]A+f_1(A,B,\chi,\tau)+
O(\mu),\\
\frac{dB}{dt}=\frac{1}{2}
\left[b(\tau)+e(\tau)-
\frac{h'(\tau)}{h(\tau)}
\right]B+f_2(A,B,\chi,\tau)+
O(\mu),
\end{array}
\label{eq15}
\end{equation}
where
$$
h(\tau)=\left[-\frac{d(\tau)}{a(\tau)}
\right]^{1/2},\quad \chi(\tau,\mu)=\frac{1}{\mu}
\int\limits_{0}^\tau[-a(s)c(s)]^{1/2}ds.
$$
The functions $f_1(A,B,\chi,\tau)$,~ 
$f_2(A,B,\chi,\tau)$ are periodic in $\chi$ with 
period $2\pi$. They contain the terms  
not lower than quadratic in $A$ and $B$.
Analysis of the system (30) gives the following
theorem for the system (14) (see similar Theorem
16.2 in [6]).
(Below, the mean value of
a periodic function $f(\tau)$ is denoted by
$\langle f(\tau)\rangle$.)

\bigskip

{\it Let a 
resonance point $J_{pq}$ meet the 
conditions of Theorem 1. Let the inequality
$$
a(\tau)d(\tau)<\sigma_1<0,\quad
t\in(-\infty,\infty).
$$
holds. Finally, let  the inequality
$$
\langle b(\tau)+e(\tau)\rangle\ne 0
$$
holds. Then system (14) in the 
$\varepsilon$-neighborhood of the resonance point, 
for sufficiently small $\varepsilon$, has a unique 
almost periodic solution. This solution is 
asymptotically stable if
$$
\langle b(\tau)+e(\tau)\rangle<0,
$$
and unstable if
$$
\langle b(\tau)+e(\tau)\rangle>0.
$$}

\bigskip

Theorems 1 and 2 apply to the study of resonant
solutions in the vibro-impact system (11).
The resonance point is defined from equation
\begin{equation}
\omega_0(J_{pq})=-\frac{\pi\Omega}{\pi-arctan\frac{J_{pq}}
{2\Omega\Delta}}=\frac{q}{p}\nu.
\label{eq16}
\end{equation}
From equation (31) follows that
$$
\omega'_0(J_{pq})>0.
$$
\par
To calculate remaining coefficients of the equation
(30) we need to specify the function 
$f(t,\tau)$.
\par
Let the function $f(t,\tau)$ is defined by the 
formula (8).
\par
To calculate $a(\tau)=f_{0\eta}(\eta_0,\tau)$
we need to average $F_0(t,\tau,\varphi,J_{pq})$. 
The first term of this function is 
\begin{equation}
-4\omega_0(J_{pq})E(\tau)\sin(\nu t+\beta(\tau))
\kappa_\psi(\varphi+\frac{q}{p}\nu t).
\label{eq19}
\end{equation}
Since
$$
\kappa_\psi(\varphi+\frac{q}{p}
\nu t)=-\omega_0(J_{pq})^{-1}
\frac{1}{\pi}
\sum_{k=1}^\infty\frac{k\sin k
(\varphi+\frac{q}{p}\nu t)}
{\Omega_0^2-k^2},
\quad \Omega_0=\Omega
[\omega_0(J_{pq}]^{-1},
$$
we have to average summands in the form
$$
E(\tau)\sin(\nu t+\beta(\tau))\sin k
(\varphi+\frac{q}{p}\nu t),
\quad k=1,2,\dots.
$$
It is easy to see that the mean 
value of function (32) will be
non-zero if and only if $q=1,\,p=n\,(n=1,2,\dots)$. 
For $q=1,\,p=n$ it equals
$$
\frac{2E(\tau)\nu^2}{\pi n(\Omega^2-
\nu^2)}\cos(n\eta-\beta(\tau)).
$$
The mean value of the second summand
$$
-4\gamma J_{pq}\omega_0^2
(J_{pq})\kappa_\psi
(\varphi+\frac{q}{p}\nu t)
\kappa_\psi(\varphi+\frac{q}{p}
\nu t)
$$
equals
$$
-\frac{\gamma J_{pq}}{2}
\left(1+\frac{4\Omega^2\Delta^2}
{J_{pq}^2}\right).
$$
In calculating the mean value, we 
used the following equality
$$
\sin^{-2}\pi\Omega_0=1+4J^{-2}
\Omega^2\Delta^2.
$$
Consequently, 
the function $\eta_0(\tau)$ is determined as 
the solution of the equation
\begin{equation}
\cos(n\eta-\beta(\tau))=\frac{\gamma J_{pq}\pi n}
{4E(\tau)\nu^2}(\Omega^2-\nu^2)\left(1+\frac{4\Omega^2
\Delta^2}{J_{pq}^2}\right)=A_n(\tau).
\label{eq20}
\end{equation}
Since  $A_n(\tau)\to\infty$ as 
$n\to\infty$, we see that equation (33)
can have solutions only for a 
finite number of the values of $n$. 
If equation (33) has solutions for 
a given value $n$ ($|A_n(\tau)|<1$), then these 
solutions are determined by the following formulas
$$
\eta_{0l}(\tau)=\frac{\beta(\tau)}{n}\pm
\frac{\arccos A_n(\tau)}{n}+\frac{2l\pi}{n},
\quad l=0,\dots,n-1.
$$
Calculating the derivative of the 
function $f_0(\eta,\tau)$ at the points
$\eta_{0l}(\tau)$ yields
\begin{equation}
f_{0\eta}(\eta_{0l}(\tau),\tau)=\pm\frac{2E(\tau)\nu^2}
{\pi(\Omega^2-\nu^2)}\sqrt{1-A_n^2(\tau)}, 
\label{eq21}
\end{equation}
and, therefore, (34) has a positive sign 
at $n$ points and a negative sign at $n$ points.
Straightforward calculations similar to the ones
done earlier show that 
$$
\langle b(\tau)+e(\tau)\rangle=\langle f_1(\eta_0,\tau)
+g_{0\eta}(\eta_0,\tau)\rangle<0.
$$
Theorem 1 and 2 imply the following result.
If the resonance point $J_{1n}$ is a solution of
equation (31), then, if $\varepsilon$ is
sufficiently small and $|A_n(\tau)|<1$, equation 
(11) has $n$ solutions in the 
$\sqrt\varepsilon$ - neighborhood of the resonance
point, which are unstable resonant almost periodic
in $t$, and $n$ solutions, which are asymptotically 
stable resonant almost periodic in $t$ in the
$\varepsilon$-neighborhood of the resonance point.
\par
The case when the force is determined by formula
(9) is more simple. The functions 
$f_0(\eta,\tau)$,~$f_1(\eta,\tau)$,~$g_0(\eta,\tau)$
do not depend on $\tau$. Thus $\eta_0(\tau)\equiv
\eta_0=const$. We obtain
$$
\cos(n\eta-\theta)=\frac{\gamma J_{pq}\pi n}
{4A\nu^2}(\Omega^2-\nu^2)\left(1+\frac{4\Omega^2
\Delta^2}{J_{pq}^2}\right)=A_n
$$
and
$$
\eta_{0l}=\frac{\theta}{n}\pm
\frac{\arccos A_n}{n}+\frac{2l\pi}{n},
\quad l=0,\dots,n-1.
$$
The variable $\tau$ is contained only in terms of 
order $O(\mu^3)$ in systems (15) and (28).
We get a result similar to result
obtained for the perturbation (8).
\section{Conclusion}
\label{x}
In this paper we have shown that in the 
considered vibro-impact system with one degree 
of freedom under small biharmonic perturbation 
with close and very different frequencies can 
occur stable and unstable almost periodic 
resonant modes.
\par
Our study is based on a combination of the method
of averaging on infinite interval [6] and 
the methods of singular perturbation theory 
[10,13]. The change of variables of method of averaging
transforms the original system into 
a system with a smooth principal part. This system is 
then singularly perturbed in the slow time. Further 
study uses some ideas from the theory of 
singular perturbations.
\par
Similar methods can be used to study the 
resonant modes arising from the action of 
biharmonic perturbation in various problems of 
theory of oscillation described by equations with 
smooth and discontinuous coefficients.

\section*{References}

1. V.S. Aslanov, Chaotic behavior 
of the biharmonic dynamics 
system,\\ International 
Journal of Mathematics and 
Mathematical Sciences, 
vol. 2009, Article
ID 319179, 18pp. 
\par
2. V.I. Babitsky, V.L. Krupenin, 
Vibrations of Strongly 
Nonlinear Discontinuous\\ 
Systems, Springer - Verlag: Berlin, Heidelberg, 
2001. 
\par
3. V.I. Babitsky, M.Z. Kolovsky,
To investigation of the resonance regimes
in vibro-impact systems, Mechanics of 
Solids, n. 4, pp. 88--91, 1976 (in Russian). 
\par
4. I.I. Blekhman, P.S. Landa, 
Conjugate resonances 
and bifurcations in nonlinear systems under
biharmonical excitation, Int. J.
Non-Lin. Mech., vol. 39(2004), 421--426. 
\par
5. V. S. Burd,  Resonant almost 
periodic oscillations 
in system with slow varying 
parameters, Int. J. Non-Lin. Mech., vol. 32, 
No. 6(1997), 1143--1152. 
\par
6. V. Burd, Method of Averaging for Differential
Equations on an Infinite 
Interval. Theory and 
Applications, Chapman\&Hall/CRC, 
Boca Raton, London,
New York, 2007.   
\par
7. V.N. Chizhevsky, Analytical study of vibrational
resonance in an overdamped bistable oscillator,
International Journal of Bifurcation and Chaos,
vol. 18, No. 6(2008), 1767--1773. 
\par
8. M. Gitterman, Bistable oscillator driven by 
two periodic fields, J. Phys. A: Math. Gen. 34(2001),
L355--L357.
\par   
9. J.A. Morrison, Resonance behavior of a 
perturbed system depending on a slow-time parameter, 
J. Math. Anal. Appl., vol. 21, No. 1(1968),
79--98. 
\par
10. R.E. O'Malley, Singular Perturbations Methods
for Ordinary Differential Equations (Applied
Mathematical Sciences, vol. 89), Springer-Verlag,
New York, 1991.
\par
11. I.G. Malkin, Some problems in the theory of
nonlinear oscillations, U.S. Atomic Energy
Comissions. Translation series, AEC-tr-3766, Oak
Ridge Tenn., 1959.
\par
12. G. P. Tolstov, Fourier Series, Prentice Hall,
Inc., Englwood Cliffs, New Jersey, 1962.
\par
13. A.B. Vasil'eva, V.F. Butuzov, Asymptotic 
Methods in Singular 
Perturbations, Vysshaya 
Shkola,
Moscow, 1991 (in Russian).
\par
14. V.F. Zhuravlev, D.M. Klimov, Applied Methods 
in the Theory Oscillations, Nauka, Moscow, 1988
(in Russian).

\end{document}